\documentclass[12pt]{amsart}
\usepackage{amsmath,amssymb,amscd}

\input xy
\xyoption{all}

\DeclareMathOperator{\ext}{Ext} 
\DeclareMathOperator{\tor}{Tor}
\DeclareMathOperator{\coker}{coker}

\DeclareMathOperator{\reg}{reg}

\DeclareMathOperator{\ann}{Ann}

\DeclareMathOperator{\ass}{Ass}

\DeclareMathOperator{\der}{Der}
\DeclareMathOperator{\Hom}{Hom}
\newtheorem{thm}{Theorem}[section]
\newtheorem{cor}[thm]{Corollary}
\newtheorem{lem}[thm]{Lemma}
\newtheorem{conj}[thm]{Conjecture}
\newtheorem{prop}[thm]{Proposition}
\theoremstyle{definition}
\newtheorem{defin}[thm]{Definition}

\newcommand{\Z}{\ensuremath{\mathbb{Z}}}

\newcommand{\A}{\ensuremath{\mathcal{A}}}
\newcommand{\D}{\ensuremath{\mathcal{D}}}
\newcommand{\om}{\ensuremath{(0:_M \mathfrak{m}^{\infty})}}

\newcommand{\p}[1]{\ensuremath{\mathbb{P}^#1}}
\newcommand{\sheaf}[1]{\ensuremath{\mathcal{#1} }}

\newtheorem{example}[thm]{Example}

\begin{document}
\title[Castelnuovo-Mumford regularity by approximation]{Castelnuovo-Mumford regularity\\by approximation}
\author{Harm Derksen}
\address{Department of Mathematics, University of Michigan, 525 East University Ave., Ann Arbor, MI 48109-1109}
\email{hderksen@umich.edu}
\thanks{The first author is partially supported by NSF, grant DMS 0102193.}
\author{Jessica Sidman}
\address{Department of Mathematics, 970 Evans Hall, UC Berkeley,  Berkeley, CA 94720}
\email{jsidman@math.berkeley.edu}
\address{Department of Mathematics and Statistics, 415 A Clapp Lab, Mount Holyoke College, South Hadley, MA 01075}
\email{jsidman@mtholyoke.edu}
\thanks{The second author is supported by an NSF Postdoctoral Fellowship, DMS 0201607}
\maketitle
\begin{abstract}
The Castelnuovo-Mumford regularity of a module gives a rough measure of its complexity.  We bound the regularity of a module given a system of approximating modules whose regularities are known.  Such approximations can arise naturally for modules constructed by inductive combinatorial means.  We apply these methods to bound the regularity of ideals constructed as 
combinations of linear ideals and the module of 
derivations of a hyperplane arrangement as well as to give degree bounds for invariants of finite groups.
\end{abstract}
\section{Introduction}\label{section:1}

Let $S = k[x_0, \ldots, x_n]$ where $k$ is an infinite field and  
let $\mathfrak{m} = (x_0, \ldots, x_n)$ be the homogeneous maximal ideal.  Throughout, all modules will be finitely generated graded $S$-modules.

Our main goal is to present a method for bounding the Castelnuovo-Mumford 
regularity of an $S$-module $M$ in terms of the 
regularities of modules that form an approximation system for $M.$  We also show how to use an approximation system to find which primes may be associated to $M.$  Precise definitions are given in \S\ref{section:3}.

For applications we use the technique of approxmation to bound the regularity of any ideal constructed by taking sums, products, and intersections of linear ideals. We also discuss how our methods lead to degree bounds for invariants of finite groups and a bound on the regularity of the module of derivations tangent to a hyperplane arrangement.

The Castelnuovo-Mumford regularity of a module $M$ governs the degrees 
appearing in a minimal graded free resolution of $M.$  Furthermore, if $M$ is in sufficiently general coordinates over a field of characteristic zero, it 
bounds the degrees needed in Gr\"obner basis computations with respect to 
the reverse lexicographic ordering on monomials (see~\cite{bayer-stillman}).
Recent efforts from several different quarters contribute to our understanding 
of how regularity behaves under operations such as addition, multiplication, 
intersection, and taking radicals (cf.~\cite{conca-herzog,chardin-dcruz}).  

The idea for our approximations is rooted in a method
of Conca and Herzog used to compute the regularity of a product of linear 
ideals (see~\cite{conca-herzog}).  Our generalization of this method recaptures the original regularity bound of Conca and Herzog as well as the bound on the regularity of intersections of linear ideals given recently by the authors in~\cite{derksen-sidman}.

In \S\ref{section:2} we provide definitions and background on regularity.
In \S\ref{section:3} we discuss  general techniques for working with approximations of modules.  We give applications to combinations of linear ideals and degree bounds for invariants of finite groups in \S\ref{section:4} and applications to modules of derivations in \S\ref{section:5}.  Some of these results also appear in the thesis of the second author in \cite{sidmanthesis}.\\[10pt]

{\bf Acknowledgements:} We would like to thank Rob Lazarsfeld, Aldo Conca, and J\"urgen Herzog for helpful conversations.  We also thank Hal Schenck for introducing us to the subject of derivations of hyperplane arrangements and Bernd Sturmfels for encouraging us to think about the primes associated to a combination of linear ideals.  Finally, we would like to acknowledge the role of the software package \emph{Macaulay 2} \cite{grayson-stillman} in computations of concrete examples as we worked on the project.  We also thank the referee for many helpful comments.

 \section{Working with Castelnuovo-Mumford regularity}\label{section:2}
In this section we review basic facts and definitions for working with regularity.  The following definition of Castelnuovo-Mumford regularity shows how the regularity of a module governs the degrees appearing in a minimal resolution:
\begin{defin}
Let $M$ be a finitely generated graded module over $S$ and let 
\[ 0 \to F_l \to \cdots \to F_0 \to M \to 0\] be a minimal graded free 
resolution of $M.$  We say that $M$ is $r$-regular in the sense of Castelnuovo 
and Mumford if $F_i$ is generated by elements of degree $\leq r+i$ for all $i$.  
The regularity of $M,$ denoted $\reg(M),$ is the least integer
$r$
 so that $M$ is $r$-regular. 
\end{defin}
Alternate formulations in terms of $\tor$, $\ext$ and local cohomology (which we will not need here) lend themselves to quick homological proofs.  In particular, one can use such methods to prove the basic fact:

\begin{lem} [Corollary 20.19 in \cite{eisenbud}]\label{lem: seq}
Let $A$, $B$, and $C$ be finitely generated graded modules over $S$ and let \[ 0
 \longrightarrow A \longrightarrow B \longrightarrow C \longrightarrow 0\] be a 
short exact sequence.  Then
\begin{itemize}
\item[(a)]  $\reg(A) \leq \max\{\reg(B), \reg(C) +1\}$.
\item[(b)]  $\reg(B) \leq \max\{\reg(A), \reg(C)\}$.
\item[(c)]  $\reg(C) \leq \max\{\reg(A) -1, \reg(B) \}$.
\end{itemize}
\end{lem}

In studying the regularity of a graded module $M$ it is often helpful to understand $\om,$ the submodule of all elements of $M$ that are killed by a power of ${\mathfrak m}.$   The reader familiar with local cohomology will note that $\om$ is just $H^0_{\mathfrak m}(M)$.  One can see from Proposition \ref{prop: hypersurface} that the regularity of $M$ is always at least as large as the regularity of $\om.$  Furthermore, since the regularity of a module of finite length is equal to its top nonzero degree, $\reg M$ is always greater than or equal to the top degree of $\om.$

The following definition provides some useful terminology for discussing 
elements that form a regular sequence modulo the finite length part of a
module.

\begin{defin}[cf.~\cite{trung}]
Let $M$ be a finitely generated graded $S$-module.  
An element $z$ is \emph{filter-regular} for $M$ if the multiplication map 
\[z: M_{i-\deg(z)} \to M_i\] is injective for all $i \gg 0.$  
A sequence $z_0, \ldots, z_m$ is a \emph{filter-regular} 
$M$-sequence if each $z_i$ is a filter-regular element for 
$M/(z_0, \ldots, z_{i-1})M.$  Note that $z$ is a filter-regular element for
$M$ if and only if it is a nonzerodivisor on 
$N = M/\om.$
\end{defin}

Proposition \ref{prop: hypersurface}, due to Conca and Herzog, describes the 
relationship between the regularity of $M$ and its quotient by a 
sufficiently general hypersurface.  It generalizes Lemma 1.8 in \cite{bayer-stillman} which is stated in terms of homogeneous ideals and a statement for modules and hyperplanes in~\cite{eisenbud}.  It is a key ingredient in the proof of 
Theorem~\ref{thm: m^k approx} which gives a bound on the regularity of a 
module that is approximated sufficiently well by modules with known 
regularities.

\begin{prop}[Proposition 1.2 in \cite{conca-herzog}]\label{prop: hypersurface}
If $x$ is a homogeneous form that is a filter-regular element for $M,$
 then \[\reg(M) =\max \{\reg \om, \reg(M/xM) -\deg(x)+1\} .\]
\end{prop}

\section{Approximations}\label{section:3}
In this section we introduce the notion of an approximation system for a module and develop techniques for using approximation systems to find associated primes and regularity bounds.

\begin{defin}\label{def: approx sys}Let $M$ be a finitely generated graded $S$-module.  We will call a finitely generated graded $S$-module $M'$ an \emph{approximation} of $M$ if there is a surjective morphism 
\[\phi: M \twoheadrightarrow M'\] with $\ker \phi$ 
annihilated by a proper homogeneous ideal $I$ of $S.$

Let $M_1, \ldots, M_d$ be finitely generated graded modules such that 
each $M_i$ is an approximation of $M$ with morphism 
$$\phi_i: M \twoheadrightarrow M_i$$
and $I_i \cdot\ker \phi_i = 0$ for homogeneous ideals $I_1, \ldots, I_d$.  
If $\sum I_i \supseteq \mathfrak{m}^t$ for some integer $t$, 
then we say the modules $M_1, \ldots, M_d$ form an \emph{approximation system} 
for $M$ of degree $t$.
\end{defin}

 In Proposition \ref{prop: assprimes} we use an approximation system to determine which primes may be associated to $M.$  Theorem \ref{thm: m^k approx} provides a first step in using approximation systems to bound $\reg(M).$  
 
\begin{prop}\label{prop: assprimes}
If $M_1, \ldots, M_d$ form an approximation system for $M$ of degree $t$ then 
\[ \ass M \subseteq \{\mathfrak{m} \} \cup \ass M_1 
\cup \cdots \cup \ass M_d.\]
\end{prop}
\begin{proof}
We begin by showing that $\ass M\cup V(I_i)=\ass M_i \cup V(I_i)$ for each $i.$

Let $f\in M$ and let $\ann(f)$ be its annihilator ideal.  If $x\in \ann(\phi_i(f))$, 
then we get $\phi_i(x f)=x \phi_i(f)=0$ and $xf\in \ker(\phi_i)$. For
any $y\in I_i$ we have $yxf=0$ so $yx\in \ann(f)$. We have shown
\begin{equation}\label{annihilators}
I_i\ann(\phi_i(f))\subseteq \ann(f)\subseteq \ann(\phi_i(f))
\end{equation}
(the right inclusion is trivial). The associated primes of $M$ and $M_i$
are the minimal primes (irreducible components) in $V(\ann(f))$ and
$V(\ann(\phi_i(f))$ respectively for various choices of $f\in M$.
Together with (\ref{annihilators}) this implies that $\ass M\cup V(I_i)=\ass M_i \cup V(I_i)$ for each $i.$

It follows that
\begin{equation}\label{assprimes}
\ass M\subseteq \ass M_i\cup V(I_i)\subseteq
\ass M_1 \cup \cdots\cup \ass M_d \cup V(I_i)
\end{equation}
for all $i$. Since ${\mathfrak m}^t\supseteq \sum_i I_i,$
we see that $\bigcap_{i=1}^d V(I_i)=\{{\mathfrak m}\}$ and
now the proposition follows from (\ref{assprimes}). 
\end{proof}
\begin{thm}\label{thm: m^k approx}
Let $M$ be a finitely generated graded module over $S$  and let 
$M_1, \ldots, M_d$ be an approximation system for $M$ of degree $t$.
Suppose that $y\in S_1$ is a linear form
 which is a filter-regular element for $M$
and on  $M_i$ for all $i$
simultaneously.
If $\reg(M_i) \leq r-t$ for all $i$ and $\reg(M/yM) \leq r-t+1$, then
we have  $\reg(M) \leq r$.
\end{thm}

\begin{proof}
By Proposition \ref{prop: hypersurface} it suffices to show that 
 $\reg \om \leq r.$  Since $\om$ has finite length it is enough to 
 show that it does not contain any elements of degree greater than $r.$  

Suppose that $f \in \om$ with $\deg(f) \ge r+1.$ 
Choose $l$ maximal so that we can write $f=y^lh$ for some $h\in M$.
Now also $h\in \om$ because $y$ is a filter-regular element for $M$.  
Let 
$$\psi_y: M \to M/yM$$
 be the quotient morphism. 
We have $\psi_y(h)\neq 0$ because of the maximality of $l$.
Since $h$ is killed by a power of ${\mathfrak m}$, so is $\psi_y(h)$
and therefore $\psi_y(h)\in (0:_{M/yM}\mathfrak{m}^\infty)$.

This shows that 
$$\deg(h)\leq \reg(M/yM)=r-t+1$$ 
and $l\geq t$.
In particular, we can write $f=y^tg$ with $\deg(g)\geq r-t+1$.
A power of ${\mathfrak m}$ kills $f=y^tg$, and therefore also
$\phi_i(y^tg)=y^t\phi_i(g)$ for each $i.$  We get that $y^t\phi_i(g)\in
(0:_{M_i}\mathfrak{m}^\infty)$ which implies that 
$\phi_i(g)\in (0:_{M_i}\mathfrak{m}^\infty)$
because of our choice of $y$.

Since
$$
\reg(0:_{M_i}{\mathfrak m}^\infty)\leq \reg M_i\leq r-t
$$
we have that $\phi_i(g)=0$.
  We obtain $g \in \bigcap_{i=1}^d \ker \phi_i.$
We can write $y^t=y_1+y_2+\cdots+y_d$ with $y_i\in I_i$ for all $i$.
Here $I_1,I_2,\dots,I_d$ are as in Definition~\ref{def: approx sys}.
We conclude that $f=y^tg=\sum_{i}y_ig=0$.
\end{proof}

Note that Theorem \ref{thm: m^k approx} holds without the assumption that the 
maps $\phi_i: M \rightarrow M_i$ are surjective.  Furthermore, the same proof goes through if we replace the assumption that we have ideals $I_i$ that annihilate the kernels of the $\phi_i$ and with the assumption that the kernel of the map $M \rightarrow \oplus M_i$ is annihilated by a power of the maximal ideal.  In certain circumstances, the weaker condition on the annihilation of the kernel of 
  $M \rightarrow \oplus M_i$ results in a better bound.  In our applications Theorem \ref{thm: m^k approx} is used in an inductive fashion which requires the stronger condition on the kernel of each $\phi_i.$

The following lemma will be helpful for the inductive application of Theorem \ref{thm: m^k approx} in what follows.

\begin{lem}\label{lem: generic 1-forms}
Let $M$ be a finitely generated graded $S$-module.  
Then there is an open dense subset $U \subseteq S_1^{n+1}$ such that each $(z_0, \ldots, z_n) \in U$ is a filter-regular $M$-sequence and $(z_0, \ldots, z_n) = \mathfrak{m}.$
\end{lem}

\begin{proof}
Let $N = M/\om.$  Since $\mathfrak{m}$ is not associated to $N,$ the set of linear forms that are filter-regular elements for $M$ is an open dense subset of $S_1.$  Here, if $N = 0$ all linear forms are filter-regular elements for  $M.$  The claim follows by induction.
\end{proof}

Notice that one must know a bound on the regularity of the quotient of a module $M$ by a sufficiently general hyperplane in order to apply Theorem \ref{thm: m^k approx}.  We eliminate this hypothesis in Theorem \ref{thm: regapprox}. 
\begin{thm}\label{thm: regapprox}
Let $M$ be a graded module which is finitely generated 
in degree $\leq r-(t-1)(n+1)$. Suppose that
$M_1,M_2,\dots,M_d$ form an approximation system for $M$ of degree $t$
such that $\reg(M_i)\leq r-(t-1)(n+1)-1$ for all $i$.
Then $\reg(M)\leq r$.

\end{thm}
\begin{proof}The result will follow from repeated applications of 
Theorem~\ref{thm: m^k approx}.  
By Lemma \ref{lem: generic 1-forms} we can choose a sequence of 1-forms 
$z_0, \ldots, z_n$ that generates $\mathfrak{m}$ and is a filter-regular sequence for
$M, M_1, \ldots, M_d$ simultaneously. 
Define $L_{j+1} = (z_0, \ldots,
z_j).$
We will prove by {\em decreasing} induction on $j$ that $\reg(M/L_{j}M)\leq 
r-(t-1)j$.

For $j=n+1$ we get $M/L_{n+1}M=M/{\mathfrak m}M$.
Since $M$ is generated in degree $\leq r-(t-1)(n+1)$,
$M/{\mathfrak m}M$ is zero in degree $> r-(t-1)(n+1)$.
It follows that
$$\reg(M/L_{n+1}M)=\reg(M/{\mathfrak m}M)\leq r-(t-1)(n+1).$$

Let us assume that $\reg(M/L_{j+1}M)\leq r-(t-1)(j+1)$. We will prove
that $\reg(M/L_{j}M)\leq r-(t-1)j$. In order to do this, we would
like to apply Theorem~\ref{thm: m^k approx}. We put $N=M/L_{j}M$
and $N_i=M_i/L_{j}M_i$ for all $i$. 
Since $L_j$ is generated by a filter-regular sequence, Proposition~\ref{prop: hypersurface} implies that
\begin{multline*}
\reg(N_i)\leq \reg(M_i)\leq r-(t-1)(n+1)-1\leq \\
\leq r-(t-1)(j+1)-1=
r-(t-1)j-t.
\end{multline*}
Also $z_j$ is a filter-regular element for $N,N_1,\dots,N_d$ and
$$
\reg(N/z_jN)=\reg(M/L_{j+1}M)\leq r-(t-1)(j+1)=r-(t-1)j-t+1.
$$
To apply Theorem~\ref{thm: m^k approx} we only need to show that
$N_1,N_2,\dots,N_d$
form an approximation system for $N$ of degree $t$.
The surjective morphism \[\phi_i:M\twoheadrightarrow M_i\] restricts to a 
surjective
morphism $L_jM\twoheadrightarrow L_jM_i$ which will be denoted by
$\widehat{\phi}_i$. Also $\phi_i$ induces a morphism 
\[\overline{\phi}_i:N=M/L_jM\to N_i=M_i/L_jM_i.\]
We have the following diagram:
\[ \xymatrix{
0 \ar[r] & L_jM \ar[r] \ar@{->>}[d]^{\widehat{\phi}_i}& M \ar[r] \ar@{->>}[d]^{\phi_i} 
& N\ar[r] \ar@{->}[d]^{\overline{\phi}_i}&0 \\
0 \ar[r] & L_jM_i \ar[r] & M_i \ar[r] & N_i \ar[r] & 0.
}\]
From the snake lemma it follows that $\overline{\phi}_i$ is also surjective and that we have an exact sequence:
$$
0\to \ker\widehat{\phi}_i\to \ker\phi_i\to \ker\overline{\phi}_i\to
0=\coker \widehat{\phi}_i.
$$
In particular, this shows that $\ker\overline{\phi}_i$ is
annihilated by $I_i$ because $\ker(\phi_i)$ is annihilated by $I_i$.
It follows that $N_1,N_2,\dots,N_d$ form an approximation system for $N$
of degree $t$.

We can apply Theorem~\ref{thm: m^k approx} and obtain
$$
\reg(N)=\reg(M/L_jM)\leq r-(t-1)j.
$$
By induction we have shown that
$$
\reg(M/L_jM)\leq r-(t-1)j.
$$
for all $j$. In particular, for $j=0$ we get
$\reg(M)\leq r$.
\end{proof}
For $t=1$, the formulation of the theorem is nicer:
\begin{cor}
Suppose that $M$ is a  graded module, finitely generated in degree
$\leq r$. Suppose that $M_1,M_2,\dots,M_d$ form an approximation system
for $M$ of degree $1$ and that $\reg(M_i)\leq r-1$ for all $i$.
Then we have $\reg(M)\leq r$.
\end{cor}
With Corollary \ref{cor: M} we show that we can eliminate the hypothesis of an upper bound on the degrees of the generators of $M$ if $M$ is a submodule of a direct sum of copies of $S.$ Note that the hypotheses of Corollary \ref{cor: M} are naturally satisfied if $M$ is a proper homogeneous ideal of $S.$

\begin{cor}\label{cor: M}
Let $F=S^l$ be a free module. Suppose that $M\subseteq F$ is a graded module and that we have graded modules $M_1,M_2,\dots,M_d\subseteq F$ and ideals $I_1,I_2,\dots,I_d$ such that
$$
I_i\cdot M_i\subseteq M\subseteq M_i\quad\forall i.
$$
If $I_1+I_2+\cdots+I_d={\mathfrak m}$, and $\reg(M_i)\leq r-1$ for some $r \ge 2$ and all $i$, then $\reg(M)\leq r$.
\end{cor}
\begin{proof}

First note that $\reg(M) \leq \max\{0, \reg(F/M)+1 \}$ from Lemma~\ref{lem: seq}.  We have a surjective morphism $\phi_i:F/M\to F/M_i$ whose
kernel is killed by $I_i$. Since $\sum_i I_i={\mathfrak m}$,
we have that $F/M_1,F/M_2,\dots,F/M_d$ form an approximation system
for $F/M$ of degree $1$. Since \[\reg(F/M_i)\leq \max\{0, \reg(M_i)-1\} \leq r-2\] and $F/M$ is generated in degree $0\leq r-1$, we get
$\reg(F/M)\leq r-1$ and $\reg(M)\leq r$.
\end{proof} 
In the applications in this paper, all approximations
of ideals and modules will be linear (i.e., of degree 1).  However, variations on the notion of an approximation of degree $t$ seem to be useful as well. Suppose that ${\mathfrak p}\subseteq S$ is a prime ideal. To understand its regularity, one might approximate the ideal ${\mathfrak p}$ with other ideals. For example, at a certain point $x\in \p{n}$, the ideal might be locally a complete intersection of $f_1,f_2,\dots,f_s$ where $s$ is the height of ${\mathfrak p}.$   The ideal $(f_1,\dots,f_s)$, being a complete intersection, will have a good regularity bound.
We can find a homogeneous $f\in S$ such that
$$f{\mathfrak p}\subseteq (f_1,\dots,f_s)\subseteq {\mathfrak p}.$$
Now $(f_1,\dots,f_s)$ could be called an approximation of ${\mathfrak p}$,
but it is not quite an approximation in the sense of 
Definition~\ref{def: approx sys}.
In fact, this type of approximation is in a sense dual to the one in
Definition~\ref{def: approx sys} and we might call it a {\it co-approximation}.
We define a co-approximation of a module $M$ as an injective morphism
$\phi:M'\to M$ such that the cokernel is killed by a proper ideal.
We can of course also define a {\it co-approximation system of degree $t$} in
a similar fashion. If ${\mathfrak p}$ is locally a complete intersection
then it will have a co-approximation system of complete intersection ideals
(which have good regularity bounds).

One can develop a theory for co-approximations that is similar to the theory for approximations presented here. We will mention only one result for co-approximations, namely Corollary~\ref{cor: coapprox} which follows from our results for approximations.  We will not pursue the theory of co-approximations here,  nor will we give regularity bounds for (locally complete intersection) prime ideals.
We just would like to point out that the scope of the methods
employed in this paper are not limited to subspace arrangements.
The method of (co)-approximation of ideals may be used to obtain regularity 
bounds of arbitrary prime ideals.
\begin{cor}\label{cor: coapprox}
Let $F=S^l$ be a free module
with graded submodules $M,M_1,M_2,\dots,M_d$.
Suppose that $I_1,I_2,\dots,I_d$ are ideals
satisfying
$$
I_i\cdot M\subseteq M_i\subseteq M
$$
for all $i$ and
$$
I_1+I_2+\cdots+I_d\supseteq {\mathfrak m}^t.
$$
If $\reg(M_i)\leq r-(t-1)n$ for all $i$, then
$\reg(M)\leq r$.
\end{cor}
\begin{proof}
We choose a basis $f_{i,1},\dots,f_{i,l_i}$ of
the degree $t$ part of $I_i$.
 Let $g$ be an arbitrary nonzero homogeneous polynomial of large degree. Define 
$$
h=g\prod_{i=1}^d\prod_{j=1}^{l_i}f_{i,j}.
$$
We have
$$
h \cdot M\subseteq (h/f_{i,j})\cdot M_i\subseteq (h/f_{i,j})\cdot M
$$
for all $i$ and $j$.
We put $M'=F/(hM)$, and $M_{i,j}'=F/(h/f_{i,j})M_i$.
The kernel of the surjective
homomorphism $M'\twoheadrightarrow M_{i,j}'$ is killed by $f_{i,j}$.
The ideal ${\mathfrak m}^t$ is generated by all $f_{i,j}$.
Therefore, the $M_{i,j}'$ form an approximation system for $M'$ of degree $t$.
Now 
\begin{multline*}
\reg(M_{i,j}') \leq \deg(h)-t+\reg(M_i)-1\leq
\deg(h)-t+(r-(t-1)n)-1\\
=(r+\deg(h)-1)-(t-1)(n+1)-1
\end{multline*}
for all $i$ and $j$.  Also note, that by taking the degree of $g$ (and therefore the degree of $h$) very large, we  get that $$(r+\deg(h)-1)-(t-1)(n+1)\geq 0.$$
Since $M'$ is generated in degree at most $0$, it will also be generated in degree at most $(r+\deg(h)-1)-(t-1)(n+1)$. From Theorem~\ref{thm: regapprox} it follows that
$$
\reg(M)=\reg(M')-\deg(h)+1\leq (r+\deg(h)-1)-\deg(h)+1=r.
$$
\end{proof}

\section{The regularity of systems of ideals generated by linear
ideals}\label{section:4}
In the course of general investigations into the regularity of products of ideals and modules, Conca and Herzog \cite{conca-herzog} discovered that for linear ideals $I_1, \ldots, I_d$, $\reg(I_1 \cdots I_d) = d.$  Using Corollary \ref{cor: M} it is actually possible to prove a very general statement about the regularity of an ideal that is a ``combination'' of linear ideals in a sense that we will make precise in the following definition.

\begin{defin}
Let $A$ be a set of linear ideals in $S$.
We define $\sheaf{C}(A)$ to be the system of ideals generated by $A$,
where
$$
\sheaf{C}(A)=\bigcup_{r=0}^\infty \sheaf{C}_r(A)
$$
with 
\begin{eqnarray*}
\sheaf{C}_0(A)&=&\{(0),S\},\\
\sheaf{C}_1(A)&=&\sheaf{C}_0(A)\cup \big\{\sum_{I\in B}I\mid B\subseteq
A\big\},
\end{eqnarray*}
and $\sheaf{C}_r(A)$ equal to:
\begin{eqnarray*}
&\{{\mathfrak a}\cap {\mathfrak b}+{\mathfrak c}\mid
{\mathfrak a}\in \sheaf{C}_a(A),{\mathfrak b}\in \sheaf{C}_b(A),
{\mathfrak c}\in \sheaf{C}_1(A),a+b=r,1\leq a,b\}\\
&\cup 
\{{\mathfrak a}\cdot {\mathfrak b}+{\mathfrak c}\mid
{\mathfrak a}\in \sheaf{C}_a(A),{\mathfrak b}\in \sheaf{C}_b(A),
{\mathfrak c}\in \sheaf{C}_1(A),a+b=r,1\leq a,b\}\\
&\cup
\{{\mathfrak a}+ {\mathfrak b}\mid
{\mathfrak a}\in \sheaf{C}_a(A),{\mathfrak b}\in 
\sheaf{C}_b(A),a+b-1=r,2\leq a,b\}\quad\mbox{for $r>1$.}
\end{eqnarray*}

Observe that for $r >1$ each $I \in \sheaf{C}_r(A)$ is constructed by taking 
sums, products, and intersections of linear ideals 
(possibly with repetition) from $A$. 

If $A$ is the set of all linear ideals, then we write
$\sheaf{C}_r=\sheaf{C}_r(A)$ and $\sheaf{C}=\sheaf{C}(A)$.%\m{I altered

\end{defin}

\begin{example}
Let $A=\{(x),(y)\}$. Then we have
$$\sheaf{C}_1(A)=\{(0),(x),(y),(x,y),S\}$$ 
and $\sheaf{C}_2(A) - \sheaf{C}_1(A)$ equal to: 
\[\{ (xy), (x^2), (x^2, y), (y^2), (y^2, x), (x^2, xy, y^2), (x^2, xy), (xy, y^2)\}.\]
\end{example}
In Theorem~\ref{thm: systems} we will prove that each ideal in 
$\sheaf{C}_r$ is $r$-regular.  Note that in particular $I_1\cdot I_2\cdots I_d$
and $I_1\cap I_2\cap \cdots\cap I_d$ lie in $\sheaf{C}_d(\{I_1,\dots,I_d\})$.
So Theorem~\ref{thm: systems} implies the main result of \cite{derksen-sidman} as well as the regularity bound on products of linear ideals in \cite{conca-herzog}.

The lemma below allows us to decompose the ideals in $\sheaf{C}_r$.
\begin{lem}\label{lem: system decomp} 
Suppose that $A=\{I_1,I_2,\dots,I_d\}$ is a set of linear ideals.
 Assume that $J \in \sheaf{C}_r(A)$ for $r \ge 1$, and $J\not\in \sheaf{C}_r(A')$ for any proper subset $A'\subset A$.
  Then there are 
   ideals $J_1, \ldots, J_d\in \sheaf{C}_{r-1}(A)$
 such that
 \begin{equation} \label{eq: linmult} 
 I_i \cdot J_i \subseteq J\subseteq J_i
 \end{equation}
 for all $i$.
  \end{lem}

\begin{proof}
The proof goes by induction on $r$. 
The intuition here is that $J$ can be expanded into a formula consisting of 
sums, products, and intersections of the linear ideals in $A$.  
We obtain $J_i$ by replacing $I_i$ with $S$ in this formula.  
The resulting ideal must contain the original ideal.  
Furthermore, since $S$ is the identity for intersection and products,
$J_i$ is a good approximation of $J$.

Suppose $r = 1$. The cases $J=(0)$ and $J=S$ are trivial.
If $J\neq (0),S$, then $J=\sum_{I\in A}I$ because of the
minimality of $A$. We can take $J_1=J_2=\cdots=J_d=S\in \sheaf{C}_0(A)$.

Suppose that $J\in \sheaf{C}_r(A)$ for $r \ge 2$. 
There are three cases:
\begin{itemize}
\item[(i)] $J={\mathfrak a}\cap {\mathfrak b}+
{\mathfrak c}$ with
${\mathfrak a}\in \sheaf{C}_a(A)$, ${\mathfrak b}\in \sheaf{C}_b(A)$,
 ${\mathfrak c}\in \sheaf{C}_1(A)$, $a+b=r$ and $1\leq a,b$
\item[(ii)]$J={\mathfrak a}\cdot {\mathfrak b}+
{\mathfrak c}$ with
${\mathfrak a}\in \sheaf{C}_a(A)$, ${\mathfrak b}\in \sheaf{C}_b(A)$,
 ${\mathfrak c}\in \sheaf{C}_1(A)$, $a+b=r$ and $1\leq a,b$
 \item[(iii)]$J={\mathfrak a}+ {\mathfrak b}$ with
${\mathfrak a}\in \sheaf{C}_a(A)$, ${\mathfrak b}\in \sheaf{C}_b(A)$,
  $a+b=r-1$ and $2\leq a,b$.
 \end{itemize}
 We will assume that we are in case~(i). The other cases go similarly.
 We can choose minimal subsets $A',B',C'\subseteq A$
 such that $\mathfrak{a}\in \sheaf{C}_a(A')$,
 $\mathfrak{b}\in \sheaf{C}_b(B'),\mathfrak{c}\in \sheaf{C}_1(C')$.
 Let $A'=\{I_1',\dots,I_l'\}$, $B=\{I_{l+1}',\dots,I_{m}'\}$
 and $C=\{I_{m+1}',\dots,I_{p}'\}$.
By minimality of $A$, we must have 
$$
\{I_1,\dots,I_d\}=A=A'\cup B'\cup C'=\{I_1',\dots,I_{p}'\}.
$$
(there may be repetition on the right-hand side). By induction we have ideals 
 \begin{eqnarray*}\mathfrak{a}_1, \ldots, \mathfrak{a}_l\in
 \sheaf{C}_{a-1}(A'),\\
 \mathfrak{b}_{l+1}, \dots,\mathfrak{b}_m\in \sheaf{C}_{b-1}(B'),\\
 \mathfrak{c}_{m+1}=\cdots=\mathfrak{c}_{p}=S\in \sheaf{C}_0(C')
\end{eqnarray*}
such that
  \[ I_i' \cdot \mathfrak{a}_i \subseteq \mathfrak{a} \subseteq 
  \mathfrak{a}_i, \quad I_i' \cdot \mathfrak{b}_i \subseteq \mathfrak{b} 
  \subseteq \mathfrak{b}_i\quad\mbox{and}\quad
  I_i'\cdot\mathfrak{c}_i\subseteq \mathfrak{c}\subseteq \mathfrak{c}_i
   \]
   for all $i$. Define
\begin{eqnarray*}
J_i&={\mathfrak a}_i\cap {\mathfrak b}+{\mathfrak c} &\mbox{for $i=1,\dots,l$,}\\
J_i&={\mathfrak a}\cap {\mathfrak b}_i+{\mathfrak c} &\mbox{for
$i=l+1,\dots,m$ and}\\
J_i&={\mathfrak a}\cap {\mathfrak b}+{\mathfrak c}_i &\mbox{for
$i=m+1,\dots,p$.}
\end{eqnarray*}
Note that $J_i\in \sheaf{C}_{r-1}(A)$ for $i=1,2,\dots,p$.
If $1\leq i\leq l$ then
$$
I_i' J_i=I_i'(\mathfrak{a}_i\cap {\mathfrak b}+{\mathfrak c})\subseteq
(I_i\mathfrak{a}_i)\cap {\mathfrak b}+{\mathfrak c}\subseteq
{\mathfrak a}\cap {\mathfrak b}+{\mathfrak c}=J\subseteq J_i.
$$
Similarly one can show for $i=l+1,\dots,p$ that $I_i'J_i\subseteq J\subseteq
J_i$.
\end{proof}
Using Lemma \ref{lem: system decomp} in combination with Corollary \ref{cor: M}, we can show:
\begin{thm}\label{thm: systems}
If $J \in \sheaf{C}_r$
then $J$ is $r$-regular.
\end{thm}
\begin{proof}
First we use induction on $n$, where $n+1$ is the number of
variables in $S=k[x_0,x_1,\dots,x_n]$.
The case $n=0$ follows immediately, because then
$$\sheaf{C}_r=\{(0),S,(x_0),(x_0^2),\dots,(x_0^r)\}$$
and it is easy to see that every element of $\sheaf{C}_r$ is $r$-regular.

We now assume $n>0$. With induction on $r$ we will prove
that every element $J \in \sheaf{C}_r$ is $r$-regular.  This is trivial if $r = 0,1,$ so we may assume $r \ge 2.$ 
Choose a minimal set $A=\{I_1,I_2,\dots,I_d\}$ such that $J\in \sheaf{C}_r(A)$.
Suppose that 
$$
I_1+\cdots+I_d\neq \mathfrak{m}.
$$
Without loss of generality (by using a linear change of coordinates)
we may asssume that
$$
I_1+\cdots+I_d=(x_0,x_1,\dots,x_s)
$$
with $s<n$. It follows that every ideal in $\sheaf{C}(A)$ is
generated by elements in the smaller ring $S':=k[x_0,x_1,\dots,x_s]$.
In other words,
$$
J=(J\cap S')\otimes_{S'}S.
$$
Since the ring extension $S\supseteq S'$
is flat, we see that $\reg(J)=\reg(J\cap S')$.
The ideal $J\cap S'$ lies in the class
$\sheaf{C}_r$ for the smaller polynomial ring $S'$.
By induction on the number of variables of the polynomial ring,
we know that $J\cap S'$ is $r$-regular, hence $J$ is $r$-regular.

Now let us assume that
$$
I_1+I_2+\cdots+I_d={\mathfrak m}.
$$
We apply Lemma~\ref{lem: system decomp} to obtain  $J_1,J_2,\dots,J_d$
with
$$
I_i\cdot J_i\subseteq J\subseteq J_i,\quad\mbox{for all $i$},
$$
and $J_1,J_2,\dots,J_d\in \sheaf{C}_{r-1}(A)$.
 By induction on $r$ we know that $\reg(J_i)\leq r-1$ for all $i.$  
We can apply Corollary~\ref{cor: M} with $F=S$, $M_i=J_i$ for all $i$
and $M=J$ to obtain $\reg(J)\leq r$.
\end{proof}

In the following proposition we show that the associated primes of any $J \in \sheaf{C}(A)$ must be linear ideals. 
\begin{prop}

Suppose $J\in \sheaf{C}(A)$ where $A=\{I_1,I_2,\dots,I_d\}$.
Then
$$
\ass(J)\subseteq \big\{\sum_{I\in B}I\mid B\subseteq A\big\}.
$$
\end{prop}
\begin{proof}
We have $J\in \sheaf{C}_r(A)$ for some $r$.
Without loss of generality we may assume that $A$ is minimal.
We will prove the proposition by induction on $r$.
If
$$I_1+I_2+\cdots+I_d\neq {\mathfrak m}$$
then we can
reduce to the polynomial ring with fewer variables.
(The proposition is easy for the polynomial ring in 1 variable.)
Therefore, let us assume that
$$
I_1+I_2+\cdots+I_d={\mathfrak m}.
$$
We can apply Lemma~\ref{lem: system decomp} to obtain
ideals $J_1,J_2,\dots,J_d\in \sheaf{C}_{r-1}(A)$
with $I_iJ_i\subseteq J\subseteq J_i$ for all $i$.
So $S/J_1,S/J_2,\dots,S/J_d$ form an approximation system
for $S/J$. By induction we know that the associated primes
of $S/J_1,\dots,S/J_d$ are contained in
$$
\big\{\sum_{I\in B}I\mid B\subseteq A\big\}.
$$
Now the proposition follows from Proposition~\ref{prop: assprimes}.
\end{proof}

In the next section we will show the connection between ideals inductively constructed from linear ideals and degree bounds for invariants of finite groups.

\subsection{Invariants of finite groups}
Suppose that $G$ is a finite group and that $V$ is an $n$-dimensional
 representation
of $G$ over the field $k$ Then $G$ also acts on the
coordinate ring $k[V]\cong k[x_1,\dots,x_n]$ and the ring of invariant polynomials
is denoted by $k[V]^G$. The ring $k[V]^G$ is finitely generated
(see~\cite{noether1,noether2}).
In fact, if we define $\beta(k[V]^G)$ as the smallest integer
$d$ such that $k[V]^G$ is generated by invariants of degree $\leq d$ then
$$
\beta(k[V]^G)\leq |G|.
$$
This was proven by Noether in \cite{noether1} if the characteristic 
of the base field $k$ is 0 or larger than the group order $|G|$.
If the characteristic of $k$ divides the group order $|G|$ then Noether's degree bound is no longer true.  For some time it was an open conjecture whether Noether's degree bound would hold in the general non-modular case, i.e., when the
characteristic of $k$ does not divide the group order $|G|$.
This was proven by Fleischmann in \cite{fleischmann} and Fogarty gave
another proof independently in \cite{fogarty}. 

Let $J\subseteq k[V]$ be the ideal generated by all homogeneous
invariants of positive degree. If $J$ is
generated in degree $\leq d$, then $\beta(k[V]^G)\leq d$ (see~\cite{derksen}).
Let us define  $B\subset V\times V$ by
$$
B=\{(v,g\cdot v)\mid v\in V,g\in G\}
$$
Let $G$ act on $V\times V$ where the action on the first
factor is trivial and the action on the second factor is as usual.
If we define the diagonal
by  $\Delta(V)=\{(v,v)\mid v\in V\}\subset V\times V$,
then
$$
B=\bigcup_{g\in G}g\cdot \Delta(V).
$$
So $B$ is the union of $|G|$ subspaces (presuming
that the action is faithful).
In \cite{derksen}, the first author made the following observation:
\begin{prop}
Let 
$${\mathfrak b}\subseteq k[V\times V]\cong
k[x_1,x_1,\dots,x_n,y_1,y_2,\dots,y_n]$$ 
be the vanishing ideal of $B$.
If ${\mathfrak b}=(f_1(x,y),f_2(x,y),\dots,f_r(x,y))$,
then $J=(f_1(x,0),\dots,f_r(x,0))$.
(Here $x=x_1,\dots,x_n$, and $ y=y_1,\dots,y_n$.)
\end{prop}
In particular, if ${\mathfrak b}$ is generated in degree $\leq d$
then so are $J$ and $k[V]^G$. This led to the following conjectures:
\begin{conj}[Subspace Conjecture, see \cite{derksen}]
If $I$ is the vanishing ideal of a union of $d$ subspaces,
then $I$ is generated in degree $\leq d$.
\end{conj}
Sturmfels made a stronger conjecture, namely:
\begin{conj}
If $I$ is the vanishing ideal of a union of $d$ subspaces,
then the Castelnuovo-Mumford regularity of $I$ is at most $d$.
\end{conj}
It was proven by the authors 
in \cite{derksen-sidman} that both conjectures are true. 
In particular, the ideal ${\mathfrak b}$ is generated in degree $\leq |G|$
and this implies that $\beta(k[V]^G)\leq |G|$ in the nonmodular case,
which gives yet another proof of the results of Fleischmann and Fogarty.  As we have already seen, Theorem \ref{thm: systems} implies Sturmfels' conjecture (and hence also the Subspace Conjecture.)

Besides the number $\beta(k[V]^G)$, there is another constant
that seems very interesting, namely the regularity $\rho_G(V)$ of the ideal $J$.
In particular, $J$ is generated in degree $\leq \rho_G(V)$ and therefore
$$
\beta(k[V]^G)\leq \rho_G(V).
$$
Since $k[V]/J$ has finite length, it is easy to see that
$\rho_G(V)$ is the smallest positive integer $d$ such that every monomial
of degree $d$ lies in $J$. From Fogarty's proof of the Noether bound
in the non-modular case it follows that
$$
\rho_G(V)\leq |G|.
$$
Again, this also follows from our results here. Let us define
$c_G(V)$ as the smallest positive integer $d$ such that ${\mathfrak b}$
lies in the class $\sheaf{C}_d$. As observed earlier, $c_G(V)\leq |G|$
because ${\mathfrak b}$ is the intersection of $|G|$ linear ideals.
Also the ideal 
$${\mathfrak b}+(y_1,\dots,y_n)=J+(y_1,\dots,y_n)$$
 lies in the class $\sheaf{C}_{c_G(V)}$. This implies that $J$ is $c_G(V)$-regular,
 so 
 $$
 \beta(k[V]^G)\leq \rho_G(V)\leq c_G(V)\leq |G|.
 $$
Schmid proved in \cite{schmid} that
Noether's bound can  only be sharp if the group $G$ is cyclic.
For a noncyclic group $G$ in 
characteristic 0, it was proven in \cite{domokos} that  $\beta_G(V)\leq
\frac{3}{4}|G|$ if $|G|$ is even and $\beta_G(V)\leq \frac{5}{8}|G|$
if $|G|$ is odd. This result was extended to arbitrary characteristic
in \cite{sezer}. 

Suppose that $G$ is abelian, with elementary divisors $d_1,d_2,\dots,d_l$.
In \cite{schmid} Schmid made the following conjecture:
\begin{conj}
$$
\beta(k[V]^G)\leq 1+\sum_{i=1}^l(d_i-1)
$$
\end{conj}
She also showed that 
$$\beta(k[V]^G)\geq 1+\sum_{i=1}^l(d_i-1)$$
for some representation $V$ so the conjecture could not be any sharper.
Schmid's conjecture is known to be true if $G$ is an abelian $p$-group
or an abelian group whose order has at most 2 distinct prime divisors
(cf~\cite{olsen1,olsen2}).

The methods in this paper lend themselves to degree bounds
for finite groups, since $\beta(k[V]^G)\leq c_G(V)$.
Bounds for $c_G(V)$ can be obtained by inspection of the
geometry of the subspace configuration $B$.
Here we will give one example where $c_G(V)$ is much
smaller then the expected bound $|G|$.

For every $g\in G$ and every subgroup $H\subseteq G$, we define
$$\Delta_{gH}(V)=\sum_{h\in gH} h\cdot \Delta(V).$$

\begin{lem}\label{lemma: intersection}
If $G$ is abelian, then
$$
\Delta_{g_1H_1}(V)\cap \cdots\cap \Delta_{g_lH_l}(V)=\{0\}
$$
if and only if for every character $\chi$ appearing in $V$
there exist $i$ and $j$ such that $\chi(H_i)=\chi(H_j)=\{1\}$
and $\chi(g_i)\neq \chi(g_j)$.
\end{lem}
\begin{proof}
Let $G$ act on $V\times V$ ($G$ acts on both factors in the
standard non-trivial way). Since $\Delta(V)$ is $G$-stable,
and $G$ is abelian, also $h\cdot\Delta(V)$ is $G$-stable.
This shows that $\Delta_{gH}(V)$ is $G$-stable for all $g\in G$
and all subgroups $H\subseteq G$. This means that
$\Delta_{gH}(V)$ is just a direct sum of its isotypic components. 
Now 
$$
\Delta_{g_1H_1}(V)\cap \cdots\cap \Delta_{g_lH_l}(V)=\{0\}
$$
if and only if
$$
\Delta_{g_1H_1}(V)\cap \cdots\cap \Delta_{g_lH_l}(V)\cap (V^\chi\times V^\chi)=\{0\}
$$
for every character $\chi$ of $G$ where $V^\chi$ is the isotypic
component of $V$ corresponding to the character $\chi$.
Now
$$
\Delta_{g_iH_i}(V)\cap (V^\chi\times V^\chi)=V^\chi\times V^\chi
$$
if $\chi(H_i)\neq \{1\}$, and 
$$
\Delta_{g_iH_i}(V)\cap (V^\chi\times V^\chi)=\{(v,\chi(g_i)v),v\in V^\chi\}
$$
if $\chi(H_i)=\{1\}$. The lemma follows.
\end{proof}
\begin{example}
Let $G=(\Z/2\Z)^n$ and let $V$ be any representation.
We prove $c_G(V)\leq n+1$ by induction on $n$.
The case $n=0$ is clear.
We may assume without loss of generality that $V$ does not have
the trivial representation as a summand.
For $g\in G$ and a subgroup $H\subset G$
we define 
$$
\Psi_{gH}(V)=\bigcup_{h\in gH}h\cdot \Delta(V).
$$
Let $H_1,H_2,\dots,H_l$ be all subgroups
of $G$ of index 2. Choose $g_1,\dots,g_l$ such that $g_i\not\in H_i$.
Now we have
$$
\Delta_{H_1}(V)\cap \cdots\cap \Delta_{H_l}(V)\cap
\Delta_{g_1H_1}(V)\cap\cdots\cap \Delta_{g_lH_l}(V)=\{0\}.
$$
by Lemma~\ref{lemma: intersection}.
We have
$$
\Psi_{gH_i}\subseteq B\subseteq 
\Psi_{gH_i}\cup \Delta_{H_i}(V).
$$
and similarly
$$
\Psi_{H_i}\subseteq B\subseteq 
\Psi_{H_i}\cup \Delta_{g_iH_i}(V).
$$
By induction, the vanishing ideals of
$\Psi_{H_i}$ and $\Psi_{gH_i}=g\cdot\Psi_{H_i}$
lie in $\sheaf{C}_{n}$, since $H_i\cong (\Z/2\Z)^{n-1}$.
Now the vanishing ideal of $B$ is linearly approximated by the vanishing
ideals of all $\Psi_{H_i}$ and all $\Psi_{gH_i}$.
This implies that $c_{G}(V)\leq n+1$.
Although $B$ is the union of $2^n$ subspaces, its regularity
is at most $n+1$.
In particular $\beta(k[V]^G)\leq n+1$.
It is not hard to prove this degree bound directly
(see also~\cite{olsen1}).
 \end{example}

\section{The regularity of derivations tangent to a hyperplane
arrangement}\label{section:5}
Arrangements of hyperplanes give rise to constructions and computations that 
are interesting from combinatorial, topological, and algebraic points of view.  In this section we will show how Corollary \ref{cor: M} gives a bound on the regularity of the module of derivations tangent to a hyperplane arrangement.  This bound is inspired by work of Schenck \cite{schenck} for line arrangements in \p{2}.  General constructions of Yuzvinsky \cite{yuzvinsky} show that the bound is optimal.

We collect together some of the basic terminology used in the literature 
and give some basic facts about the module of derivations associated to a 
hyperplane arrangement in \S\ref{5.1} for the readers not familiar with the subject.  
The bound on the regularity of the module of derivations is proved in
\S\ref{5.2}.
\subsection{Hyperplane arrangements: basic facts and definitions}\label{5.1}
Let $S = k[x_0, \ldots, x_n]$ and let $\mathfrak{m} = (x_0, \ldots, x_n).$  For the convenience of the reader we collect together some facts and 
definitions that are well known.  An expanded version of the section 
appears in \cite{sidmanthesis}.  For a complete introduction to the subject, 
see \cite{orlik-terao}.

\begin{defin}
We define a \emph{hyperplane arrangement} in \p{n} to be a reduced subscheme 
consisting of $d$ distinct hyperplanes, $H_1, \ldots, H_d$.  We will denote such an arrangement by $\sheaf{A}$.  We will let $f_1, \ldots, f_d$ denote a choice of linear equations with $f_i$ cutting out $H_i$ and let $F = f_1 \cdots f_d$ be the generator of the ideal of the arrangement.  (Note that an arrangement of hyperplanes in $\p{n}_k$ is equivalent to a \emph{central} hyperplane arrangement in a vector space of dimension $n+1$ over $k$, i.e., a hyperplane arrangement in which each hyperplane contains the origin.)  
\end{defin}

\begin{defin}  A hyperplane arrangement $\A$ is 
\emph{linearly general} if every subset 
of $n+1$ of the defining linear forms is linearly independent.  We say that $\A$ is \emph{essential} if the defining linear forms span $\mathfrak{m}$.
\end{defin}

Recall that the {\em module of derivations}, $\der(S) \subseteq \Hom_k(S,S),$ is
the set of all $\theta\in \Hom_k(S,S)$ 
satisfying the Leibniz rule: 
\[ \theta(fg) = f\theta(g) + g \theta(f) \quad \forall f,g \in S.\] 

For example, $\der(S)$ is the free module of rank $n+1$ 
generated by the partial derivatives $\frac{\partial}{\partial x_i}$:
\[ \der(S) = S \cdot \frac{\partial}{\partial x_0} \oplus \cdots \oplus S \cdot \frac{\partial}{\partial x_n}.\]
We may identify $\der(S)$ with $S^{n+1}$ and this equips
$\der(S)$ with a grading.

We will focus our attention on the following module in the next section: 

\begin{defin}  The module of derivations tangent to a hyperplane arrangement 
$\sheaf{A}$ with defining polynomial $F$ is denoted by 
$\sheaf{D}(\sheaf{A}) = \sheaf{D}$ and is defined: $\D := \{ \theta \in \der(S)| \theta(F) \in (F)\}.$
\end{defin}

\subsection{The regularity of the module of derivations}\label{5.2}

Independently, Rose and Terao (cf.~\cite{rose-terao}) and 
Yuzvinsky (cf.~\cite{yuzvinsky}) have constructed minimal free 
resolutions of $\D$ in the case where the 1-forms defining the 
hyperplanes in the arrangement are linearly general.  It follows from 
Yuzvinsky's construction that for a linearly general arrangement of $d$
 hyperplanes in $\p{n}$, $\reg(\D) = d-n,$ and there exist arrangements 
 of $d$ hyperplanes in $\p{n}$ for which $\reg(\D) = d-1$. 
  Indeed, the \emph{boolean arrangement} defined by the ideal $(x_0 \cdots x_n)$ is an arrangement with $\reg(\sheaf{D}) =d-n = 1.$ 
 It is also  easy to construct examples of hyperplane arrangements with regularity equal to $d-1$:

\begin{example}
Let \sheaf{A} be an arrangement of points in \p{1} with defining equation $F \in S = k[x_0, x_1]$ of degree $d$.  Then $F$ also lives in $S' = k[x_0, \ldots, x_n]$ and defines an arrangement $\sheaf{A}'$ of $d$ hyperplanes in \p{n}.  (This arrangement is just a cone over $d$ points in \p{1}.)  Then \[ \sheaf{D}(\sheaf{A}') = (\sheaf{D}(\sheaf{A})\otimes_S S') \oplus S' \cdot \frac{\partial}{\partial x_2} \oplus \cdots \oplus S' \cdot \frac{\partial}{\partial x_n}.\]
\end{example}
\noindent Schenck has shown that for lines in $\p{2}$, $\D$ is $(d-2)$-regular under the assumption that the arrangement is essential \cite{schenck}.  It would be interesting to see if this result extends to higher dimensions.  A naive application of our methods of proof fail to recover this bound because the condition of being essential is not stable under the deletion of arbitrary hyperplanes.

In this section we show that for arbitrary arrangements of $d$ hyperplanes in $\p{n}$ the module $\D$ is $(d-1)$-regular.  The result follows as a corollary to Corollary \ref{cor: M}

\begin{thm}  Let $\A$ be an arrangement of $d\geq 2$ hyperplanes in 
$\p{n}$ with defining equation $F = f_1 \cdots f_d$.  
Then $\D$ is $(d-1)$-regular.

\end{thm}

\begin{proof}  
We prove the result with induction on $n$, where
$n+1$ is the number of variables in the polynomial
ring $S=k[x_0,x_1,\dots,x_n]$. For $n=1$, we have
$$
\D=S\cdot \left(\frac{\partial F}{\partial x_1}
\frac{\partial}{\partial x_0}-
\frac{\partial F}{\partial x_0}
\frac{\partial}{\partial x_1}\right)\oplus
S\cdot\left(x_0\frac{\partial}{\partial x_0}+x_1\frac{\partial}{\partial
x_1}\right)
$$
and we see that $\D$ is $(d-1)$-regular.

Suppose that
$$(f_1,f_2,\dots,f_d)\neq {\mathfrak m}.
$$
Without loss of generality we may assume that
$$
(f_1,f_2,\dots,f_d)=(x_0,x_1,\dots,x_s)
$$
with $s<n$. Put $S'=k[x_0,x_1,\dots,x_s]$
and let $\D'$ be the module of derivations along $F=0$
in $\p{s}$. We have
$$
\D=(\D'\otimes_{S'}S)\oplus 
S\cdot \frac{\partial}{\partial x_{s+1}}\oplus \cdots
\oplus S\cdot \frac{\partial}{\partial x_{n}}.
$$
By induction on $n$ we already know that $\D'$ is $(d-1)$-regular,
and therefore $\D$ is $(d-1)$-regular.

We prove now by induction on $d$ that $\D$ is $(d-1)$-regular.
If $d<n$, then $(f_1,f_2,\dots,f_d)\neq {\mathfrak m}$ and
we are done.
Suppose now that $d\geq n$ and $(f_1,f_2,\dots,f_d)={\mathfrak m}$.
Let $\D_i = \D( \A -{H_i}).$ 
These modules will play the role of the modules $M_i$ in the statement of 
Corollary \ref{cor: M}. Note that the case $d = 2 = n$ is trivial because then we may assume the hyperplanes are defined by the vanishing of coordinates.  By induction on $d$, we may assume that 
$\D_i$ is $(d-2)$-regular.   
We claim that 
\[ f_i \cdot \D_i \subseteq \D \subseteq \D_i.\]  
In fact, let $g_i = f_1 \cdots \hat{f_i} \cdots f_d$. 
 If $\theta \in \D$, then 
 \[ \theta(f_i\cdot g_i) = f_i \cdot \theta(g_i) + g_i \cdot \theta (f_i) 
 \in (f_i \cdot g_i) \subseteq (g_i).\] 
  Since $g_i$ divides one term in the sum, it must divide both. 
   But if $g_i$ divides $f_i \cdot \theta(g_i)$, since we are in a unique 
   factorization domain, $g_i$ must divide $\theta(g_i)$.  
   Therefore, $\theta \in \D_i$.  
   For the other inclusion, suppose that $\phi \in \D_i$.  
   Then \[ f_i \cdot \phi(f_i \cdot g_i) =  f_i^2 \cdot \phi(g_i) + 
   f_i \cdot g_i \cdot \phi (f_i)\] is clearly in $(f_i \cdot g_i)$.  
   Therefore, $f_i \cdot \D_i \subseteq \D$.
   We can apply Corollary~\ref{cor: M} and obtain
   $\reg(\D)\leq d-1$.
   \end{proof}

\bibliographystyle{amsplain}

\end{document}